\pdfoutput=1
\documentclass{ifacconf}

\usepackage{graphicx}      
\usepackage{natbib}        
\usepackage{amsmath}
\usepackage{amssymb}
\usepackage{enumitem}
\usepackage{algorithm}
\usepackage{algpseudocode}
\usepackage{tikz}
\usepackage{pgfplots}
\pgfplotsset{compat=1.17}
\usepackage{subfigure}

\usepackage{graphics} 
\usepackage{epsfig} 
\usepackage{epstopdf}
\usepackage{mathptmx} 
\usepackage{times} 
\usepackage{amsmath} 
\usepackage{amssymb}  
\usepackage{subfigure}
\usepackage[utf8]{inputenc}
\usepackage[english]{babel}
\usepackage{dsfont}
\usepackage{color}

\usepackage{multicol}

\newtheorem{ex}{Example}

\newcommand{\mcl}[1]{\mathcal{#1}}

\newcommand{\R}{\mathbb{R}}
\newcommand{\x}{\mathbf{x}}

\renewcommand{\S}{\mathbb{S}}
\newcommand{\N}{\mathbb{N}}

\DeclareMathOperator{\sign}{sign}

\begin{document}
\begin{frontmatter}

\title{Model Predictive Bang-Bang Controller Synthesis via Approximate Value Functions} 


\author[First]{Morgan Jones} 
\author[First]{Yuanbo Nie} 
\author[Second]{Matthew M. Peet}

\address[First]{ Department of Automatic Control and Systems Engineering,	The University of Sheffield (e-mail: morgan.jones@sheffield.ac.uk \& y.nie@sheffield.ac.uk).}
\address[Second]{School for the Engineering of Matter, Transport and Energy, Arizona State University, Tempe, AZ, 85298 USA (e-mail: mpeet@asu.edu)}

\begin{abstract}             
 {:} In this paper, we propose a novel method for addressing Optimal Control Problems (OCPs) with input-affine dynamics and cost functions. This approach adopts a Model Predictive Control (MPC) strategy, wherein a controller is synthesized to handle an approximated OCP within a finite time horizon. Upon reaching this horizon, the controller is re-calibrated to tackle another approximation of the OCP, with the approximation updated based on the final state and time information. To tackle each OCP instance, all non-polynomial terms are Taylor-expanded about the current time and state and the resulting Hamilton-Jacobi-Bellman (HJB) PDE is solved via Sum-of-Squares (SOS) programming, providing us with an approximate polynomial value function that can be used to synthesize a bang-bang controller.
\end{abstract}


\end{frontmatter}

\section{Introduction}
A bang-bang controller exclusively operates at the extremities of admissible inputs, toggling between upper and lower bounds. This type of controller is used in numerous physical systems with only binary actuation capabilities. Consider, for instance, a thermostat which can only alternate between ``on" and ``off" states. For more general systems where inputs are constrained to be within some rectangular set, bang-bang controllers are the optimal solution to a large class of Optimal Control Problems (OCPs) that have input affine dynamics and costs. Because of their optimality, Bang-bang controllers are ubiquitous in a wide domain of applications from medicine~\citep{ledzewicz2002optimal}, semiconductor gas discharge~\citep{kim2001high}, spacecraft maneuvers~\citep{taheri2018generic}, etc. Over the years many methods to synthesize bang-bang controllers have been proposed including~\cite{jacobson1970computation,dadebo1998computation,kaya2004computations}

Numerical methods for solving OCPs can be broken down into two distinct categories, direct methods and indirect methods. Direct methods parameterize the system state and control trajectories by a finite sum of basis functions on a time discretisation mesh to convert the continuous-time infinite-dimensional OCP into a Nonlinear Programming (NLP) problem that can be solved numerically. There are many toolboxes using the direct approach including ICLOCS~\citep{nie2018iclocs2} and GPOPs~\citep{patterson2014gpops}, with the subsequent NLPs solved using solvers such as SNOPT~\citep{gill2005snopt} and IPOPT~\citep{biegler2009large}. Direct methods have demonstrated impressive capabilities. However, as noted in~\cite{aghaee2021switch}, when it comes to OCPs that have bang-bang solutions direct methods may struggle due to ``ill conditioning and discontinuities in the optimal control at the switching points". Such computational issues of direct methods are also discussed in~\cite{pager2022method}. Another obstacle associated with the direct method lies in the difficulty of finding the global optimal solution for the resulting NLP. Mitigating the risk of converging to local minima entails initializing the algorithm with a good initial guess at the solution, also known as warm starting. 

Indirect methods for solving OCPs can be further broken down into two sub-categorizes, those based on Pontryagin Maximum Principle (PMP) and those based on Dynamic Programming (DP). The PMP provides necessary conditions for the optimality of an open-loop controller while DP is able to provide necessary and sufficient conditions for the optimality of a closed-loop controller through the Hamilton Jacobi Bellman (HJB) PDE. A more detailed discussion of the methods can be found in~\cite{liberzon2011calculus}. The PMP method involves minimizing the Hamiltonian along the solution map of the adjoint equation. If the Hamiltonian is convex (implying a unique optimal controller), it has been shown in~\cite{preininger2018convergence} that bang-bang controllers can be synthesized using the PMP. However, in general, this assumption is restrictive.

For these reasons, we focus on the DP approach to solving the optimal control problem. The disadvantage of using DP in continuous time, of course, is that it requires us to solve the HJB PDE -- a notoriously difficult nonlinear PDE. Various methods for approximately solving the HJB PDE exist such as the relatively recent approaches of~\cite{kalise2018polynomial} or~\cite{garcke2017suboptimal}. These methods typically require some sort of discretization of state and time. Alternatively, it is possible to approximately solve the HJB PDE by relaxing the equation to an inequality and applying the moment method~\citep{zhao2017control,korda2016controller,kamoutsi2017infinite}, Sum-of-Squares (SOS) programming~\citep{ichihara2009optimal,jennawasin2011performance} or other approximation schemes~\citep{sassano2012dynamic}. The advantage of these HJB relaxation approaches is that you can certify properties of the resulting approximate solution to the HJB PDE, such as being a uniformly lower bound. For this reason, we adopt the Sum-of-Squares (SOS) approach to computing HJB sub-value functions as described in~\cite{jones2020polynomial}. This approach  yields an approximate polynomial value function. A key advantage of this method is that by increasing the degree of the SOS problem, this function can be made arbitrary close (under the $L^1$-norm) to the true solution of the HJB PDE. 

Unfortunately the work of~\cite{jones2020polynomial} requires the vector field to be polynomial, which is not true for many applications. The contribution of this work, then, is to extend~\cite{jones2020polynomial} to tackle non-polynomial OCPs. Specifically, we propose to Taylor expand non-polynomial terms and solve the resulting approximate OCP. Unfortunately, of course, accuracy of the Taylor expansion can only be guaranteed in a small neighborhood of the expansion point and degrades quickly as we move away from this point. To address this issue, we propose a Model Predictive Control (MPC) type framework where after a small time horizon, the Taylor expansion is re-computed and the controller is re-synthesized. \\
\textbf{Notation:} For $x \in \R^n$ we denote $||x||_2= \sqrt{ \sum_{i =1}^n x_i^2 }$. Let $C^1(\Omega, \Theta)$ be the set of continuous functions with domain $\Omega \subset \R^n$ and image $\Theta \subset \R^m$. We denote the space of polynomials $p: \R^n \to \R$ by $\R[x]$ and polynomials with degree at most $d \in \N$ by $\R_d[x]$. We say $p \in \R_{2d}[x]$ is Sum-of-Squares (SOS) if there exist $p_i \in \R_{d}[x]$ such that $p(x) = \sum_{i=1}^{k} (p_i(x))^2$. We denote $\sum_{SOS}^d$ to be the set of SOS polynomials of at most degree $d \in \N$ and the set of all SOS polynomials as $\sum_{SOS}$. \vspace{-0.02cm}
\section{Solving OCPs Using Dynamic Programming} \vspace{-0.12cm}
 Consider the following family of Optimal Control Problems (OCPs), each initialized by $(x_0,t_0) \in \R^n\times [0,T]$,
\begin{align}  \nonumber 
&V(x_0,t_0):=\inf_{ u, x}   \int_{t_0}^T c(x(t),  u(t) ,t) dt + g(x(T))  \\ \label{opt: optimal control probelm}
 &\text{ subject to, }  \dot{x}(t) = f(x(t),  u(t)) \quad \text{ for all } t \in [t_0,T],\\ \nonumber
	&   ( u(t),x(t)) \in U \times \Omega\quad  \text{ for all } t \in [t_0,T], \quad   x(t_0)=x_0,
\end{align}
where $c: \R^n \times \R^m \times \R \to \R$ is referred to as the running cost; $g: \R^n \to \R$ is the terminal cost; $f: \R^n \times \R^m \to \R^n$ is the vector field; $\Omega \subset \R^n$ is the state constraint, $U \subset \R^m$ is the input constraint set; and $T$ is the final time. For simplicity, throughout this paper, we assume the existence of unique solutions to the Ordinary Differential Equation (ODE) described by the nonlinear dynamics $\dot{x}=f(x(t),u(t))$, for every admissible input to the Optimal Control Problem (OCP) and initial condition $x_0 \in \Omega$. This assumption is non-restrictive, as various well-established conditions ensure the existence and uniqueness of solution maps. For instance, when the vector field, $f$, is Lipschitz continuous, the solution map exists over a finite time interval. Moreover, this interval can be extended arbitrarily if the solution map remains within a compact set, as elaborated in~\cite{khalil2002nonlinear}.

Solving OCPs is a formidable challenge due to their lack of analytical solutions. One approach to solve OCPs, often referred to as the Dynamic Programming (DP) method,  reduces the problem to solving a nonlinear PDE. More specifically, rather than solving Opt.~\eqref{opt: optimal control probelm} we solve the following nonlinear PDE:
	\begin{align}  \nonumber
	&\nabla_t V(x,t) + \inf_{u \in U} \left\{ c(x,u,t) + \nabla_x V(x,t)^T f(x,u) \right\} = 0 \\
	& \hspace{3.75cm} \text{ for all } (x,t) \in \R^n \times [t_0,T] \label{eqn: general HJB PDE}\\ \nonumber
	& V(x,T)= g(x) \quad \text{ for all } x \in \R^n.
\end{align}
Upon solving the HJB PDE~\eqref{eqn: general HJB PDE} we can derive a controller as shown in the following theorem.

\begin{thm}[\cite{liberzon2011calculus}] \label{thm: value functions construct optimal controllers}
	Consider the family of OCPs in Eq.~\eqref{opt: optimal control probelm} with $\Omega=\R^n$ (no state constraints). Suppose $V \in C^1(\R^n \times \R, \R)$ solves the HJB PDE \eqref{eqn: general HJB PDE}. Then $ u^*: [t_0,T] \to U$ solves the OCP initialized at $(x_0,t_0) \in \R^n \times [0,T]$ with associated dynamical trajectory $x^*:[t_0,T] \to \R^n$ if and only if
\begin{align} \label{u}
		& u ^*(t)= k(x^*(t), t), 
  \text{ }  \dot{x^*(t)}=f(x^*(t),u^*(t)) \text{ for } t \in [t_0,T],\\ \label{eqn:k}
		& \text{where }	k(x,t) \in \arg \inf_{u \in U} \{ c(x,u,t) + \nabla_x V(x,t)^T f(x,u) \}.
	\end{align}
\end{thm}
Note that Theorem~\ref{thm: value functions construct optimal controllers} requires the solution to the HJB PDE to be differentiable. In practice this condition is often relaxed using a generalized notion of a solution to the HJB PDE, called a viscosity solution~\citep{bardi1997viscosity}. Also note that in the state unconstrained case, $\Omega=\R^n$, the solution to the HJB PDE corresponds to the optimal objective function of the OCP and is referred to as the value function. 

\underline{\textbf{Bang-Bang control:}} Given a solution to the HJB PDE~\eqref{eqn: general HJB PDE}, we can construct an optimal state feedback controller according to Eq.~\eqref{eqn:k}. However, Eq.~\eqref{eqn:k} is an optimization problem in itself requiring us to compute the infimum over $u \in U$. Solving this optimization at every time-step during implementation could be impractical. Fortunately, for OCPs with input-affine dynamics and costs, this auxiliary optimization problem has an analytical solution.

Consider OCP~\eqref{opt: optimal control probelm}, where the cost function is of the form $c(x,u,t)= c_0(x,t) + \sum_{i=1}^m c_i(x,t)u_i$, the dynamics are of the form $f(x,u)=f_0(x) + \sum_{i=1}^m f_i(x)u_i$, there are no state constraints, $\Omega=\R^n$, and the input constraints are of the form $U=[a_1,b_1] \times ... \times [a_m,b_m]$. WLOG we assume $a_i=-1$ and $b_i=1$ for $i \in \{1,\dots,m\}$, that is $U=[-1,1]^m$, since we can always make the coordinate substitution $\tilde{u}_i = \frac{2u_i - 2b_i}{b_i - a_i}$ for $i \in \{1,...,m\}$. Substituting this into Eq.~\eqref{eqn:k} we obtain 
{ \begin{align} \label{optimal controller 3}
		&	k(x,t) \in \arg \hspace{-0.2cm} \inf_{u \in [-1,1]^m}\bigg\{ \sum_{i=1}^m c_i(x,t)u_i + \nabla_x V(x,t)^T f_i(x)u_i \bigg\}.
			\end{align} }

The objective function in Eq.~\eqref{optimal controller 3} is linear in the decision variables $u \in \R^m$, and since the constraints have the form $u_i \in [1,-1]$, it follows that Eq.~ \eqref{optimal controller 3} can be analytically solved,
 \begin{align} \label{optimal controller new 2}
		&k_i(x,t)= - \sign(c_i(x,t)+ \nabla_x V(x,t)^T f_i(x)).
\end{align}

\section{Approximate Solutions of the HJB PDE} \vspace{-0.2cm}
Let us view the problem of solving the HJB PDE~\eqref{eqn: general HJB PDE} through the lens of optimization theory and consider this problem as a feasibility optimization problem with two equality constraints: the HJB PDE itself and the boundary condition. In optimization theory, when faced with a challenging non-convex problem, a common tactic is to relax the constraints of the optimization problem to make the feasible set convex. In the context of solving the HJB, the equality constraints make the feasible set non-convex. However, we may relax the equality constraints to inequality constraints in the following manner,
	\begin{align} \label{feas: convex  find VF}
		& \text{Find } J \in C^1(\R^n \times \R, \R) \text{ subject to:} \\
		\label{ineq: diss ineq for sub sol of HJB}
	& \nabla_t J(x,t) + c(x,u,t) + \nabla_x J(x,t)^T f(x,u) \ge 0 \\ \nonumber
	& \hspace{3cm} \text{for all } (x,u,t) \in \Omega \times U \times (t_0,T), \\ \label{ineq: BC}
	& J(x,T) \le g(x) \text{ for all } x \in\Omega .
\end{align}
Now, if $V$ is a solution to the HJB PDE then it satisfies Eqs~\eqref{ineq: diss ineq for sub sol of HJB} and~\eqref{ineq: BC} since 
\begin{align*}
0&=\nabla_t V(x,t) + \inf_{u \in U} \left\{ c(x,u,t) + \nabla_x V(x,t)^T f(x,u) \right\}\\
& \le \nabla_t V(x,t) +  c(x,u,t) + \nabla_x V(x,t)^T f(x,u) \\
& \hspace{3cm} \text{for all } (x,u,t) \in \Omega \times U \times (t_0,T).
\end{align*}
Problem~\eqref{feas: convex  find VF} is linear in the decision variable $J$. However, a function $J$, feasible for Problem~\eqref{feas: convex  find VF}, may be arbitrarily far from the value function. For instance, in the case $c(x,u,t) \ge 0$ and  $0 \le g(x)  <M$, the constant function $J(x,t)\equiv - C$ is feasible for any $C>M$. Thus, by selecting sufficiently large enough $C>M$, we can make $||J -V||$ arbitrary large, regardless of the chosen norm,  $|| \cdot ||$. To address this issue, we propose a modification of Problem~\eqref{feas: convex  find VF}, wherein we include an objective of minimizing the $L^1$ distance to the true solution to the HJB PDE~$ \int_{\Lambda \times [t_0,T]}  |V(x,t)-J(x,t)| dx dt$, where $V$ is the unknown solution to the HJB PDE $\Lambda \subset \R^n$ is some compact set. 
	\begin{align} \label{opt: convex  find VF}
	& \inf_{J \in C^1(\R^n \times \R, \R)} \int_{\Lambda \times [t_0,T]}  |V(x,t)-J(x,t)| dx dt  \text{ subject to:} \\ \nonumber
	& \nabla_t J(x,t) + c(x,u,t) + \nabla_x J(x,t)^T f(x,u) \ge 0 \\ \nonumber
	& \hspace{3cm} \text{for all } (x,u,t) \in \Omega \times U \times (t_0,T), \\ \nonumber 
	& J(x,T) \le g(x) \text{ for all } x \in\Omega .
\end{align}

Unfortunately, since $V$ is unknown we cannot simply solve Opt.~\eqref{opt: convex  find VF}. Fortunately, as we will see in the next proposition, feasible solutions to Problem~\eqref{opt: convex  find VF} are ``sub-value functions" -- i.e. functions that uniformly lower bound the true value function.  Then $\int_{\Lambda \times [t_0,T]}  |V(x,t)-J(x,t)| dx dt= \int_{\Lambda \times [t_0,T]}  V(x,t) dx dt -\int_{\Lambda \times [t_0,T]}  J(x,t) dx dt$ and since $\int_{\Lambda \times [t_0,T]}  V(x,t) dx dt$ is a constant it can be eliminated from the objective function.

\begin{prop} \label{prop: diss ineq implies lower soln}
Suppose $J \in C^1(\R^n \times \R, \R)$ satisfies Eqs~\eqref{ineq: diss ineq for sub sol of HJB} and~\eqref{ineq: BC} and $\Omega$ is compact. Then
	\begin{align*}
	&	J(x,t) \le V(x,t) \text{ for all } (x,t) \in \Omega \times [t_0,T],
	\end{align*}
	 where $V$ is given by the objective function of the OCP~\eqref{opt: optimal control probelm}.
\end{prop}
\begin{pf} 
Let $(x_0,t_0) \in \Omega \times [0,T]$. First suppose there is no feasible input to the OCP, then $V^*(x_0,t_0) = \infty$. Clearly in this case $J(x_0,t_0) < V^*(x_0,t_0)$ as $J$ is continuous and therefore is finite over the compact region ${\Omega} \times [0,T]$. Alternatively if there exists a feasible input, ${u}$, let us denote the resulting solution map of the underlying dynamics of the OCP by $\tilde{x}$, where $\tilde{x}(t) \in \Omega$ for all $t \in [t_0,T]$ and $\tilde{x}(t_0)=x_0$. By Eq.~\eqref{ineq: diss ineq for sub sol of HJB} we have $\text{for all } t \in [t_0,T]$
	\begin{align*}
	& \nabla_t J(\tilde{x}(t),t) + c(\tilde{x}(t), {\tilde{u}} (t),t) + \nabla_x J(\tilde{x}(t),t)^T f(\tilde{x}(t), {\tilde{u}}(t) ) \ge 0.
	\end{align*}
Now, using the chain rule we deduce
		\begin{align*}
	& \frac{d}{dt} J(\tilde{x}(t),t) + c(\tilde{x}(t), {\tilde{u}} (t),t)  \ge 0 \text{ for all } t \in [t_0,T].
	\end{align*}
	Then, integrating over $t \in [t_0,T]$, and since $J(\tilde{x}(T),T) \le g(\tilde{x}(T))$ by Eq.~\eqref{ineq: BC}, we have
	\begin{align} \label{above ineq}
	J(x_0,t_0) \le \int_{t_0}^{T} & c(\tilde{x}(t),{\tilde{u}}(t),t) dt + g(\tilde{x}(T)).
	\end{align}
	Since Eq.~\eqref{above ineq} holds for any feasible input, we may take the infimum over all feasible inputs to show that $J(x_0,t_0) \le V(x_0,t_0)$.
\end{pf}

For the case that $c$, $g$ and $f$ are polynomials, $U=[-1,1]^m$ and $\Omega=\{x \in \R^n: R^2 -||x||_2^2 \ge 0\}$, and $h_\Omega$ are polynomials, we are now able to eliminate $V$ from Opt.~\eqref{opt: convex  find VF} and tighten the optimization problem to the following Sum-of-Squares (SOS) optimization problem,
\begin{align} \label{opt: SOS for sub soln of finite time}
	& P_d \in \arg \max_{P \in \R_d[x]} c_{\text{f}}^T \alpha\\ \nonumber
	& \text{subject to: } k_{0},k_{1} \in \sum_{SOS}^d, \quad s_{i} \in {\sum_{SOS}^d} \text{ for } i=0,1,\dots,m+2 \\ \nonumber
	& P (x,t) = c^T_f Z_d(x,t), k_{0}(x) = g(x)-{P(x,T)} - s_{0}(x)(R^2 \hspace{-0.6pt} - \hspace{-0.6pt} ||x||_2^2),\\ \nonumber
	& k_{1}(x,u,t) =  \nabla_t P(x,t) +c(x,u,t) + \nabla_x P(x,t)^T f(x,u)    \\ \nonumber
	& \quad  -s_{1}(x,u,t)(R^2  - ||x||_2^2)  - s_{2}(x,u,t) (t-t_0)(T-t)\\ \nonumber
 & \quad -\sum_{i=1}^m s_{2+i}(x,u,t)(u_i+1)(1-u_i),
\end{align}
where $\alpha_i=\int_{\Lambda \times [t_0,T]}  Z_{d,i}(x,t) dx dt$, $Z_d: \R^n \times \R \to \R^{ {n+1 +d \choose d}}$ is the vector of monomials of degree $d \in \N$ and $c_f \in \R^{ {n+1 +d \choose d}}$ is the monomial coefficient vector. 

We next show that by solving Opt.~\eqref{opt: SOS for sub soln of finite time} we can approximate value functions, solutions to the HJB PDE~\eqref{eqn: general HJB PDE}, to arbitrary accuracy with respect to the $L^1$ norm.

\begin{prop} \label{prop: SOS converges}
	Consider OCP~\eqref{opt: optimal control probelm}, where $\Omega=\R^n$ (no state constraints), $c$, $g$ and $f$ are polynomials, $U=[-1,1]^m$ and $T>0$. If this OCP has a solution $(x^*(t),u^*(t))$ such that for some $R>0$ we have that
 \begin{align} \label{eq: condition unconstrained-constrained}
    & ||x^*(t)||_2^2 < R^2 \text{ for all } t \in [0,T],\\
 \label{convergence in L1 for SOS}
	&	\text{ then } \lim_{d \to \infty} \int_{\Lambda \times [t_0,T]}  |V(x,t) -P_d(x,t)| dx dt  = 0,
	\end{align}
	where $V$ solves the HJB PDE~\eqref{eqn: general HJB PDE} and $P_d$ is the solution to Opt.~\eqref{opt: SOS for sub soln of finite time} for $d \in \N$.
\end{prop}
\begin{pf}
    The condition given in Eq.~\eqref{eq: condition unconstrained-constrained} ensures the value function of the state constrained ($\Omega=\{x\in \R^n: ||x||_2 <R\}$) and unconstrained versions ($\Omega=\R^n$) of the OCP are equivalent. Then, following the mollification arguments of~\cite{jones2020polynomial} we approximate the Lipschitz continuous value function of the unconstrained OCP by a smooth function while satisfying the inequalities given Eqs~\eqref{ineq: diss ineq for sub sol of HJB} and~\eqref{ineq: BC}. We then approximate this smooth function by a polynomial using the Weierstrass approximation theorem. We show using Putinar's Positivstellensatz~\citep{putinar1993positive} that this polynomial is feasible to Opt.~\eqref{opt: SOS for sub soln of finite time}. Hence, by optimality, the solution to Opt.~\eqref{opt: SOS for sub soln of finite time} is closer to the value function under the $L^1$ than this feasible polynomial and this feasible polynomial can be made arbitrarily close to the value function, thus showing convergence. 
\end{pf}

\vspace{-0.1cm}
\section{MPC For Solving Non-Polynomial OCPs} \vspace{-0.15cm}
We have seen that we can approximately solve the HJB PDE~\eqref{eqn: general HJB PDE} by numerically solving Opt.~\eqref{opt: SOS for sub soln of finite time} for OCPs with polynomial costs and dynamics. However, many practical problems are not polynomial, limiting the broader applicability of the method. To rectify this issue, faced with a non-polynomial OCP, we propose to first approximate the non-polynomial terms using a Taylor expansion about the initial condition of the OCP -- i.e. $(x_0,t_0)$. More generally, we denote the Taylor operator as 
\begin{align*}
&	\mcl{T}_{y,d}f(x)\\
 &:= \sum_{k_n=0}^d \cdots \sum_{k_1=0}^d \frac{(x_1-y)^{k_1}}{k_1!} \cdots \frac{(x_n-y)^{k_n}}{k_n!}  \left(\frac{\partial^{k_1+\dots k_n} f}{\partial x_1^{k_1} \dots \partial x_n^{k_n}} \right)(y).
\end{align*}
Now, given an OCP~\eqref{opt: optimal control probelm}, with input affine dynamics and cost, we denote the polynomial terms with a superscript ``$p$" and non-polynomial terms with an ``$np$" superscript as follows,
\begin{align} \nonumber
c(x,u,t)&:=c_0^p(x,t) + \sum_{i=1}^m c_i^p(x,t)u_i + c_0^{np}(x,t) +  \sum_{i=1}^m c_i^{np}(x,t)u_i.\\ \label{eq: nonpoly cost and dynamics}
g(x)&:=g^p(x)+ g^{np}(x). \\ \nonumber
f(x,u)&:=f_0^p(x) + \sum_{i=1}^m f_i^p(x)u_i + f_0^{np}(x) +  \sum_{i=1}^m f_i^{np}(x)u_i.
\end{align}
We then replace non-polynomial terms with their degree $k$ Taylor expansions about the initial state and time conditions of the OCP ($(x_0,t_0)$), that is replace $c_i^{np}(x,t)$ with $\mcl{T}_{[x_0,t_0],k}c_i^{np}(x,t)$  and similarly for terms $f^{np}_i$ and $g^{np}$. This yields an OCP parameterized exclusively by polynomials and which can be solved using the methods described in the preceding section -- yielding a bang-bang feedback control law.

Unfortunately, the resulting value function only matches the true value function near the initial condition of the OCP, $(x_0,t_0)$. Because the feedback law is based on the minimization of this value function, we expect that the performance of the resulting feedback controller will match the true optimal controller near the initial condition, but will diverge as the solution evolves. Therefore, after implementing the controller over an implementation period of length $T_I>0$ we propose to re-synthesize the controller based on a new approximated OCP with non-polynomial terms expanded about a new initial condition $(x(t_0+T_I),t_0+T_I)$ -- yielding a new approximate value function that closely matches the true value function about the Taylor expansion point taken further along the systems trajectory and hence synthesising a new approximately optimal feedback law. Fig.~\ref{fig:MPC} illustrates how we expect the Taylor expansion error to vary over time. 

Finally, in order to reduce complexity of solving the SOS optimization problem, and inspired by the MPC framework, we take a receding horizon approach and only synthesize the controller over a reduced prediction horizon length $T_h\ge T_I$ (See Fig.~\ref{fig:MPC}). Of course, unless $T_I=T_h=T$, if the receding horizon solution is to match the solution to the original OCP, we must necessarily assume $T=\infty$ and in this case, we require $g(x)=0$ -- i.e., there is no terminal cost. Note, however, that since we allow for time-varying dynamics and cost functions, we do not eliminate time as a dependent variable in the HJB or value function. Moreover, because these prediction horizons may be small we may wish to only approximate the value function over some reduced integration region that depends on the final state of the previous implementation period, that is $\Lambda=[\delta_1(x),\delta_2(x)]^n$ where $\delta_i: \R^n \to \R$. The resulting algorithm is illustrated in Fig.~\ref{fig:MPC} and summarized in Algorithm~\ref{alg}.
\begin{algorithm} 
\caption{Polynomial VF-based Bang-Bang MPC}
 \textbf{Input:}\\
	 \text{OCP parameters: }$c,f$ of Form~\eqref{eq: nonpoly cost and dynamics}, $T>0, x_0 \in \R^n$   \\
 \text{MPC parameters:} $T_h>0$, $T_I>0$. \\
 \text{SOS parameters:} $d \in \N$, $\delta_1$, $\delta_2$ $R>0$.\\
 \text{Taylor parameter:} $k \in \N$. 
	\begin{algorithmic}[1]
		\For{$t \in \{0,T_I,2T_I, \dots, T-T_I\}$}
		\State $\tilde{c}(x,u,t)=c_0^p(x,t) + \sum_{i=1}^m c_i^p(x,t)u_i + \mcl{T}_{[x_0,t],k}c_0^{np}(x,t) +  \sum_{i=1}^m \mcl{T}_{[x_0,t],k}c_i^{np}(x,t)u_i  $ \Comment{Replace running cost by Taylor approx cost}
						\State $\tilde{f}(x,u)=f_0^p(x) + \sum_{i=1}^m f_i^p(x)u_i + \mcl{T}_{x_0,k}f_0^{np}(x) +  \sum_{i=1}^m \mcl{T}_{x_0,k}f_i^{np}(x)u_i  $ \Comment{Replace dynamics by Taylor approx}
		\State Compute $P_d$ by solving Opt.~\eqref{opt: SOS for sub soln of finite time} for $\tilde{c},\tilde{f},R>0$ and $\tilde{g}(x) \equiv 0$ over $[t,t+T_h]$ with integration region $\Lambda=[\delta_1(x_0),\delta_2(x_0)]^n$.
		\State $k_i(x,t)= - \sign(c_i(x,t)+ \nabla_x P_d(x,t)^T f_i(x))$ \Comment{Eq.~\eqref{optimal controller new 2}}
		\State \text{Simulate } $\dot{x}^*(t)=f(x^*(t),k(x^*(t),t))$ over $t \in [t+T_I]$.
		\State $x_0=x^*(T_I)$ \Comment{Set initial condition to be terminal trajectory}
		\EndFor
	\end{algorithmic}  \label{alg}
\end{algorithm}
	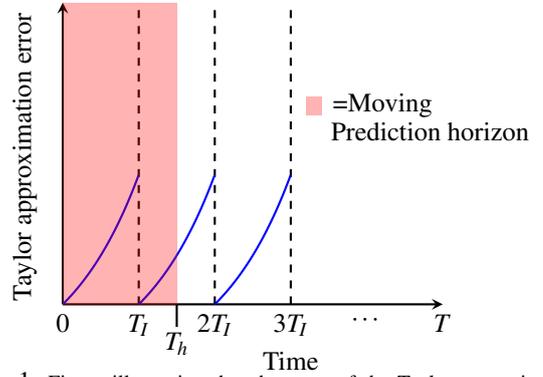
\begin{figure}
	\centering
	\begin{tikzpicture}[>=stealth, thick]
	
	\draw[->] (0,0) -- (5,0);
	\draw[->] (0,0) -- (0,4);
	
	\draw[dashed, black, thick] (2,0) -- (2,4);
	\draw[dashed, black, thick] (1,0) -- (1,4);
	\draw[dashed, black, thick] (3,0) -- (3,4);
	
	\foreach \x/\label in {0/$0$, 1/$T_I$, 2/$2T_I$, 3/$3T_I$, 4/$\cdots$, 5/$T$}
	\node[below] at (\x,0) {\label};
	
	\draw[black, thick] (1.5,0) -- (1.5,-0.3);
	\node[below] at (1.5,-0.25) {$T_h$};		
	\draw[blue, domain=0:1, samples=100] plot (\x,{(exp(\x)-1});
	\draw[blue, domain=1:2, samples=100] plot (\x,{(exp(\x-1)-1});
	\draw[blue, domain=2:3, samples=100] plot (\x,{(exp(\x-2)-1});

	\node[left, rotate=90] at (-0.5,4) {Taylor approximation error};
	\node[below] at (3,-0.5) {Time};
	
	\fill[red, opacity=0.3] (0,0) rectangle (1.5,4);
	
\fill[red, opacity=0.3] (3.2,2.5) rectangle (3.4,2.75) node[right, opacity=1, align=center] {\\ \\  \textcolor{black}{\hspace{-1.25cm}=Moving} \\ \textcolor{black}{Prediction horizon}};

	
\end{tikzpicture}
\vspace{-15pt}
		\caption{{\footnotesize Figure illustrating that the error of the Taylor approximation in Alg.~\ref{alg} grows with time, resetting after the implementation time, $T_I$, has passed. The controller is designed over a prediction horizon of length $T_h$, which could exceed the implementation time $T_I$. }} \label{fig:MPC}
			
\end{figure}

\vspace{-0.15cm}
\section{Numerical Examples} \vspace{-0.15cm}
We next present several numerical examples of using Algorithm~\ref{alg} to solve OCPs. To evaluate the performance  we approximate the objective/cost function over each implementation period using the Riemann sum:
\begin{align} \label{R sum}
\int_0^T c(\phi_{{f}}(x_0,t,&  u),t) dt   \approx  \Delta t\sum_{i=1}^{T/T_I}\sum_{j=1}^{N-1}  c(x_i(t),{u}_i(t_{i,j}),t_{i,j}) ,
\end{align} 
where $(i-1)T_I=t_{i,1}<...<t_{i,N}=iT_I$, $\Delta t = t_{i,j+1}-t_{i,j}$ for all $i \in \{1,...,N-1\}$, ${u}_i$ is the controller synthesized for implementation over $[(i-1)T_I,iT_I]$, and $x_i(t)$ can be found using Matlab's \texttt{ode45} function to solve $\dot{x}(t)=f(x(t),{u}_i(t))$.

All SOS programs are solved by utilizing Yalmip~\citep{lofberg2004yalmip} in conjunction with the SDP solver Mosek~\citep{aps2019mosek}. It's worth noting that in every numerical example, we scale the dynamics to confine the state within the range of $[-1,1]^n$. This scaling strategy ensures that both the objective function and constraints of the associated SDP problem remain within relatively small values. Consequently, this prevents the SDP solver from terminating prematurely due to suspected unboundedness. 

We benchmark Algorithm~\ref{alg} against ICLOCS~\citep{nie2018iclocs2} where we consider two numerical examples with non-polynomial time varying costs and dynamics derived from practical systems, both of which have a two dimensional state space and a one dimensional input space. However, it should be noted that both methods are not limited to systems of these dimensions and can be applied to arbitrary dimensions. Both methods are fundamentally different, being of the indirect and direct class of solution methods, but numerical examples demonstrate competitive performance. Both methods have multiple parameters that can be used to improve performance at the expense of computation time, so the results are not necessarily indicative of overall performance. Specifically, ICLOCS requires an initial solution guess whereas Algorithm~\ref{alg} is based on solving a sequence of convex optimization problems and therefore requires no such initialization. Therefore, Algorithm~\ref{alg} could be used to help warm start direct methods like ICLOCS.  

\begin{figure} 
	\centering
 \subfigure[\footnotesize Phase plot for Example~\ref{ex:van poll}.]{\includegraphics[width=0.48 \linewidth,trim = {1cm 1cm 1cm 1cm}, clip]{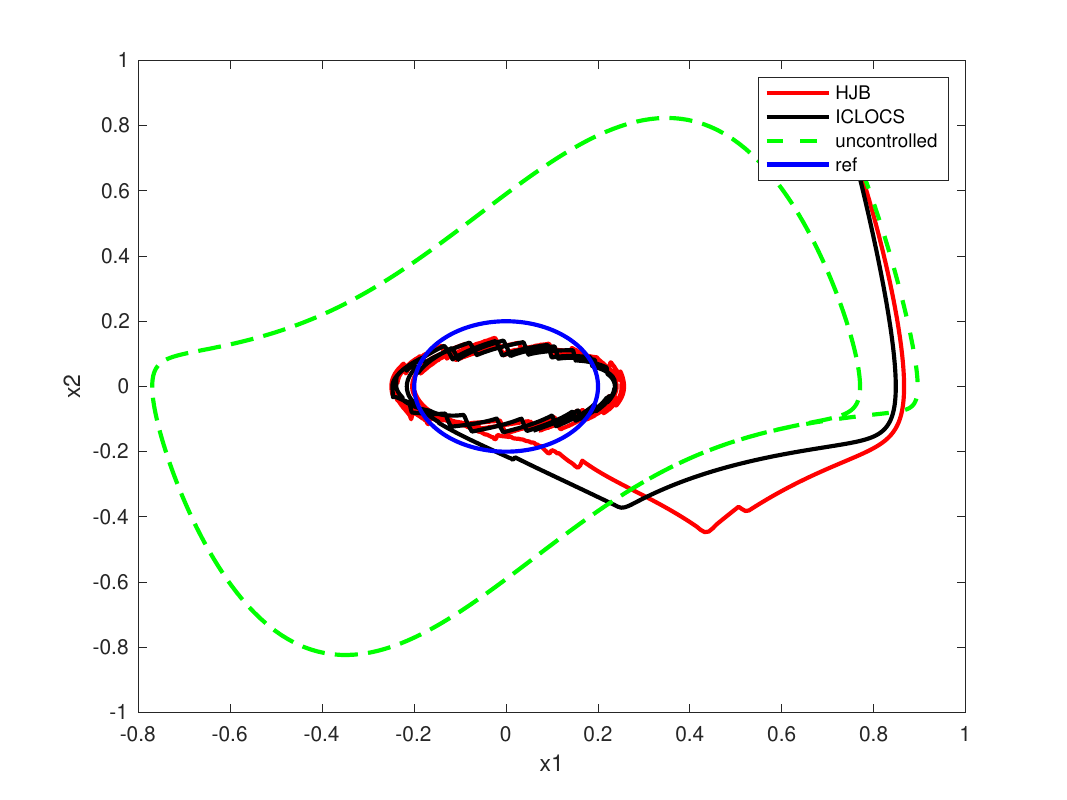} \label{fig: phase van}} 
  \subfigure[\footnotesize Phase plot for Example~\ref{ex: smib}.]{\includegraphics[width=0.48\linewidth,trim = {1.5cm 1cm 1cm 1cm}, clip]{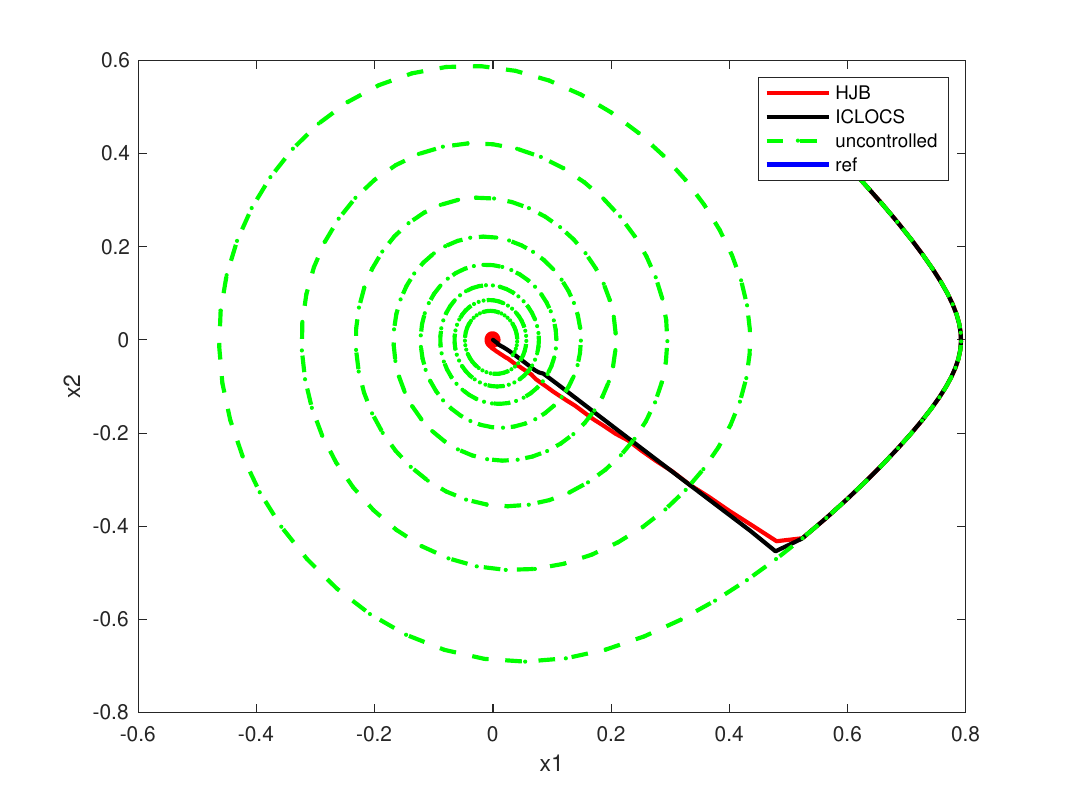} \label{fig: phase smib} } 
 \vspace{-10.5pt}
	\caption{ \footnotesize Phase plots for Examples~\ref{ex:van poll} and~\ref{ex: smib}.}  \vspace{-0pt}
\end{figure}

\begin{figure} 
\centering 	\includegraphics[width=0.9 \linewidth,trim = {0.9cm 1cm 1cm 1.25cm}, clip]{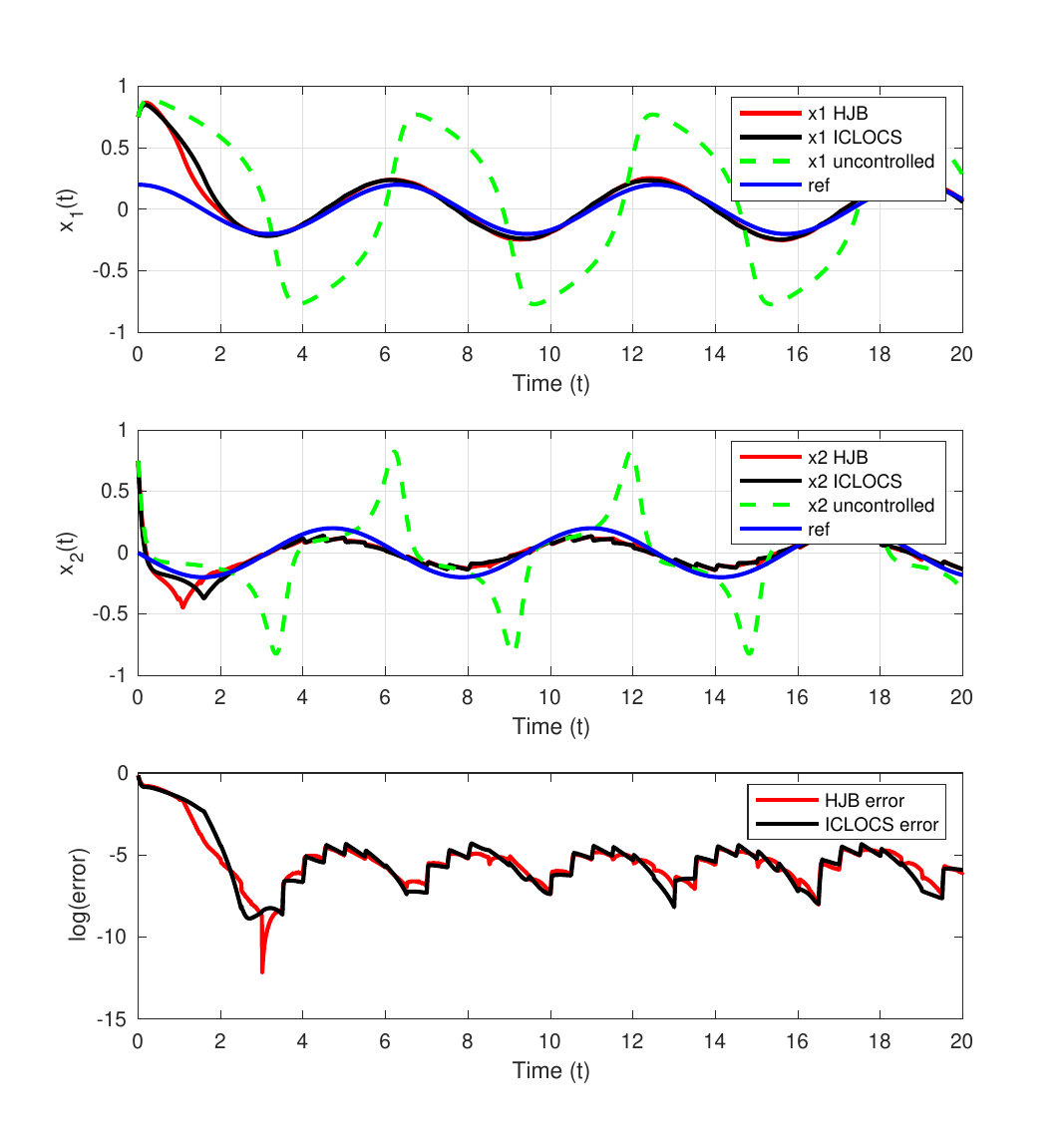}
\vspace{-12pt}		\caption{\footnotesize Figure associated with Example~\ref{ex:van poll} showing how the controller effects individual states.} \label{fig: ind van}
\end{figure}

\begin{ex}[Controlling the Van der Pol oscillator] \label{ex:van poll} \text{ }\\
Consider the following OCP
	\vspace{-0.05cm}\begin{align} \label{OCP: van poll}
		& \inf_{{u}} \bigg\{\int_0^T ||x(t)-y_{ref}(t)||_2^2dt 
  \bigg\}\\\nonumber
		&\text{Subject to: } u(t) \in[-1,1], \quad x(0)=[0.75, 0.75]^\top \\ \nonumber
		& \begin{bmatrix}
				\dot{x}_1(t) \\ \dot{x}_2(t) 
		\end{bmatrix}= \begin{bmatrix}
		2x_2(t) \\ 10x_2(t)(0.21-1.2^2x_1(t)^2) -0.8x_1(t)
	\end{bmatrix}		+ \begin{bmatrix} 0 \\ 1 \end{bmatrix} u(t)\\ \nonumber
& y_{ref}(t) = \begin{bmatrix} 0.2\cos(-t) \\0.2 \sin(-t) \end{bmatrix}.
	\end{align}
The goal of OCP~\eqref{OCP: van poll} is to minimize the tracking error between the state and some parametric clockwise circular curve of radius $0.2$. Although the dynamics of OCP~\eqref{OCP: van poll} are polynomial the cost function is not. 

The terminal time is set as $T=20$, controller implementation time is set to $T_I=0.5$ and prediction horizon is set to $T_h=1$. We use Alg.~\ref{alg} to solve this OCP with the SOS program degree set to $d=5$, the Taylor expansion degree set to $k=4$, the integration region set to $\Lambda=[\delta_1,\delta_2]^n$, where $\delta_1=-\delta_2=-0.75$, and the computation domain given by $R=||[0.75,0.75]^\top||_2$. Performance was analyzed using Eq.~\eqref{R sum} with simulation time-step $\Delta t=0.01$ and a cost of $0.521206$ was found. In contrast, for the same $T_I$ and $T_h$, $51$ discretization points, and degree $3$ polynomial approximations of state trajectories and inputs, ICLOCS incurred a slightly larger cost of $0.560883$. Fig.~\ref{fig: phase van} shows the phase space plot and the slight difference in the trajectories resulting from the two different methods. How the individual states vary with respect to time is given in Fig.~\ref{fig: ind van} as we as a log plot of how the running cost, $||x(t)-y_{ref}(t)||_2^2$, varies with respect to time.
\end{ex}


\begin{figure} 
\centering
	\includegraphics[width=0.9 \linewidth,trim = {0.9cm 1cm 1cm 1.25cm}, clip]{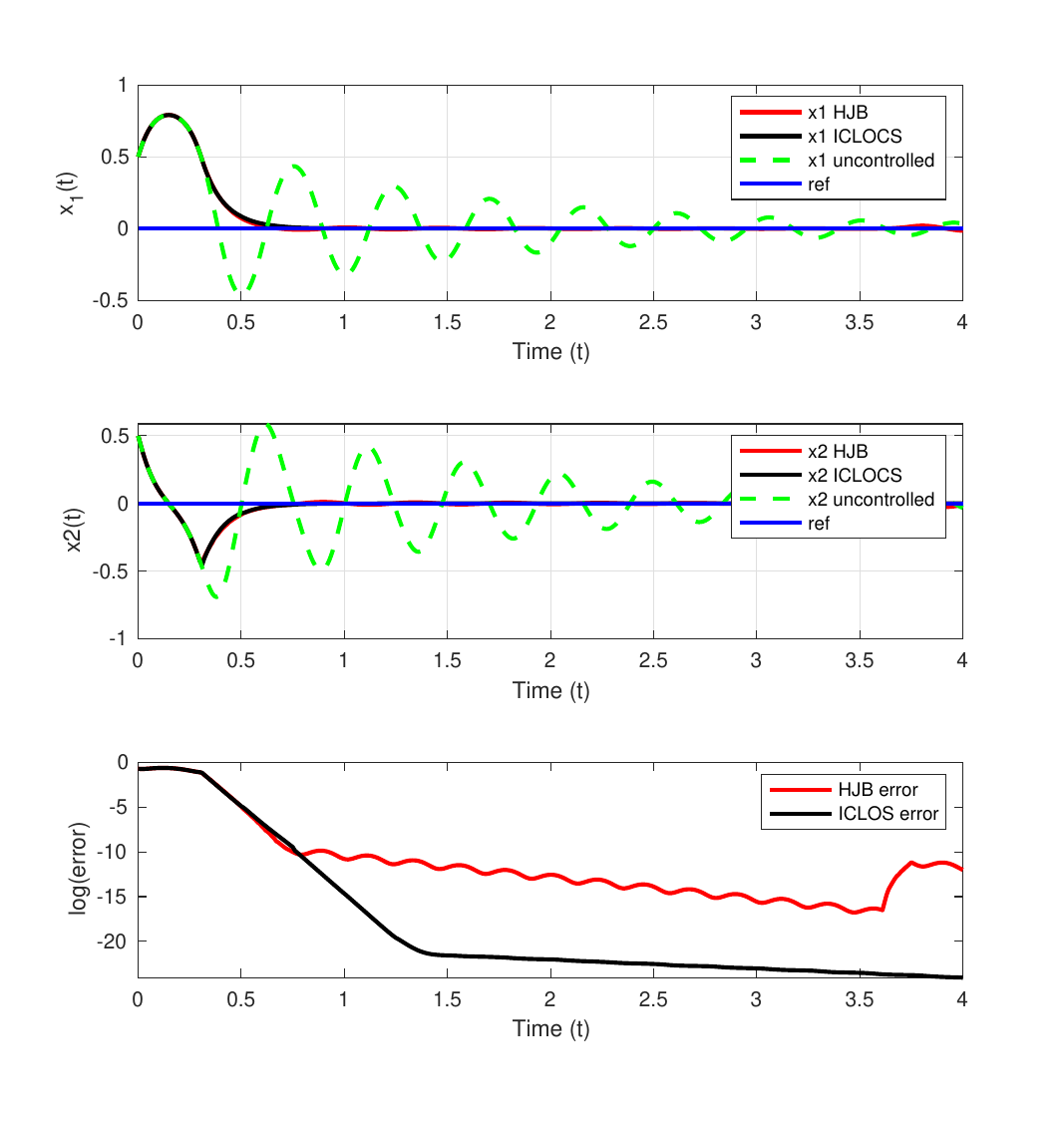} \vspace{-17pt}		\caption{ \footnotesize Figure associated with Example~\ref{ex: smib} showing how the controller effects individual states.} \label{fig: ind smib} \vspace{-8pt}	
\end{figure}

\begin{ex}[The Single Machine Infinite Bus (SMIB)] \label{ex: smib} \text{ }\\
Improving the transient stability of the SMIB system  has previously been studied in~\cite{chang1998time,ford2006nonlinear} using the PMP. We also solve this problem by considering the following OCP,
\vspace{-0.05cm}	\begin{align} \label{OCP: SMIb}
				& \inf_{{u}} \bigg\{\int_0^T e^{-t} ||x(t)||_2^2dt  \bigg\}\\\nonumber
						&\text{Subject to: } u(t) \in[-1,1], \quad x(0)=[1.5, 15]^\top \\ \nonumber
		&\begin{bmatrix}	\dot{x}_1(t) \\ 	\dot{x}_2(t) \end{bmatrix} =  \begin{bmatrix}
			x_2(t)\\ \frac{P_m-Dx_2(t)}{2H}
		\end{bmatrix} -  \begin{bmatrix}
		0\\ \frac{P_e}{2H}\sin(x_1(t)+\delta_{ep})
	\end{bmatrix}u(t),
	\end{align}
	where $H=0.0106$, $X_t=0.28$, $P_m=1$, $E_s=1.21$, $V=1$, $P_e=(E_s V)/(P_m X_t)$, $D=0.03$ and $\delta_{ep}=\sin^{-1}(1/P_e)$.
	
	Before solving the OCP we scale the dynamics to evolve over $[-1,1]^2$ by making the coordinate transformation $x \to Lx$ where $L=\begin{bmatrix}
		3 & 0 \\ 0 & 30
	\end{bmatrix}$. The goal of this problem is to improve the transient stability -- i.e. rate of convergence to the origin.  Note that the problem has non-polynomial dynamics and cost.
 
 The terminal time is set as $T=4$, controller implementation time is set to $T_I=0.25$ and prediction horizon is set to $T_h=0.5$. We use Alg.~\ref{alg} to solve this OCP with the SOS program degree set to $d=6$, the Taylor expansion degree set to $k=5$, the integration region set to $\Lambda=[\delta_1(x),\delta_2(x)]^n$, where $\delta_1(x)=x-0.2$ and $\delta_1(x)=x+0.2$, and the computation domain given by $R=1$. Performance was analyzed using Eq.~\eqref{R sum} with simulation time-step $\Delta t=0.01$ and a cost of $0.1672$ was found. In contrast, for the same $T_I$ and $T_h$, $51$ discretization points, and degree $3$ polynomial approximations of state trajectories and inputs, ICLOCS incurred a slightly smaller cost of $0.1671$. Fig.~\ref{fig: phase smib} shows the phase space plot and the slight difference in the trajectories resulting from the two different methods. How the individual states vary with respect to time is given in Fig.~\ref{fig: ind smib} as we as a log plot of how the running cost, $e^{-t}||x(t)||_2^2$, varies with respect to time.
 
\end{ex}

\vspace{-0.4cm}
\section{Conclusion} \vspace{-0.3cm} In this paper we have presented an iterative method for solving optimal control problems based on sequentially the approximated problem, whose dynamics and costs are Taylor expanded about the terminal state and time of the previous iteration. During each iteration, the HJB PDE is approximately solved using convex optimization and the resulting approximated value function is used to synthesize a bang-bang controller. Numerical examples demonstrate that this approach is competitive with the current state of the art direct method solver. Unlike direct methods, our method does not require an initial guess of the solution. Future work includes the possibility of precomputing a library of approximated value functions for different meshes of the state space and time interval offline. Then rather than sequentially solving the HJB PDE in an online MPC fashion the controller can be rapidly synthesized based on the current state and time by selecting the appropriate value function from this library.
\vspace{-0.5cm}

\bibliography{ifacconf.bib} 

                                                   







\end{document}